\documentclass{amsart}

\usepackage{amssymb}
\usepackage{latexsym, amsmath,amsthm}

\theoremstyle{plain}

\theoremstyle{definition}

\def\Z{\mathbb{Z}}

\begin{document}

\title {Branching Data for Curves up to genus 48}

\author{Jennifer Paulhus} 
\date{\today }
\address{Department of Mathematics and Statistics, Grinnell College, 
Grinnell, IA 50112}
\email{paulhusj@grinnell.edu}
\thanks{Partially supported by a Harris Faculty Fellowship through Grinnell College}


\begin{abstract} An algorithm of Thomas Breuer produces complete lists of
  automorphism groups of curves for a fixed genus, and the code to execute
  this algorithm was written in the computer algebra package GAP and run by
  Breuer over a decade ago. For each curve $X$ of genus $g$ and a
  group $G$ acting on $X$,   the branching data for the covering $X \to X/G$ is
  computed but was not record.
  We have added functionality to the code to output this data, and have
  translated the code to the computer algebra package Magma. This article
  provides an explanation of the code to add this branching
data. Complete data from this modified program for curves up to genus 20,
and data for genus up to $48$ of groups $G$ so that $|G|>4(g-1)$ are
available.  We also include sample programs which
provide a way to search the data for groups or actions with  desired
properties, these functions should be of use to researchers who want to find examples
of group actions on lower genus curves with certain properties, or who need
explicit examples of generating vectors for their research.  The code also
may be used to create generating vectors for higher genus curves, if the
automorphism group and signature of the mapping $X \to X/G$ are already
know.  
 \end{abstract}

\maketitle


\section{Introduction}\label{S:Intro}

Thomas Breuer devised an algorithm to generate a list of all groups acting
on a Riemann Surface of a given genus \cite{breuer}. The algorithm depends on databases
of groups of a given order in a computer algebra package such as Magma
\cite{magma} or GAP \cite{gap}.  In the late 1990's, Breuer coded
the algorithm in GAP and ran it for genus up to
48.  His code and the data have been available upon request.

The list of groups includes the full automorphism group of any
curve of that genus, as well as all possible subgroups.  Work of Magaard
et al. \cite{g4-10monodromy} and Bujalance et al. \cite{bujalance}
determines which groups are the full automorphism groups for genera up to $10$ and all hyperelliptic
curves, respectively. Additional work on automorphism groups of
hyperelliptic curves may be found in \cite{brandtstich} and
\cite{shaska}. Conder has used a routine in Magma to generate so called
``large'' automorphism groups for genera up to $101$ \cite{conder}.

In the process of computing the groups, Breuer's
algorithm also determines the branching data of a covering, as well as the
conjugacy classes in which the ramification occurs for that covering. But
this information was not recorded
in the original data set. This additional branching data is vital for a technique for decomposing Jacobian
varieties \cite{paulhus}, and this was the author's initial motivation for 
extending Breuer's code to also produce the branching data.

Additionally, Breuer outlines an algorithm for
expanding the results to groups of order beyond the largest order present
in small group databases from computer algebra packages.  The key distinction between this algorithm
and the lower genus one is that this algorithm does not need to search
through all groups of a fixed order. Instead, it requires a database of
perfect groups, and recursively uses the lower genus data.  This means
having searchable files containing all the data for lower genus will
be important for implementing the higher genus algorithm.

Such an expanded database could be useful in projects where a result about
automorphism groups could be reduced to low genus cases, and then proven
by an exhaustive search through the expanded database (see
\cite{woottonfullautogp} or \cite{woottonprimecovers}, for example). The extended
database could also be used to answer a variety of polynomial
problems involving Galois extensions (for instance,
generalizing a result of Fried \cite[Proposition 2]{fried}). 
 
\vspace{.1 in}

We have translated Breuer's code into Magma, and modified it to also record
the generating vectors of the coverings. For a separate project, we have begun to
use this modified code to implement the higher genus algorithm. An outline
of Breuer's algorithm  is
provided in Section \ref{S:Background}.  In Section \ref{S:code} we describe the Magma code used to generate the
data discussed above, and we explain some slight modifications made to
Breuer's original data files. Descriptions of  the files attached to this paper may be found in Section
\ref{S:data}, and some additional functionality is described in Section
\ref{S:other}, particularly how to use the provided programs to search
through the data.  All files described in this paper may be found here:

\verb+http://www.math.grinnell.edu/+$\sim$\verb+paulhusj/monodromy.html+

\section{Background}\label{S:Background}

 Below is a
summary of  the main ideas of the algorithms in Breuer's work.  More
details may be found in \cite{breuer}.  \\ 

The  automorphism group of a Riemann surface is the group consisting of all bijective holomorphic maps from the surface to itself.  This group has long been known to be finite and bounded in size by  $84(g-1)$ where $g$ is the genus of the surface.  

Let $\mathbb{H}$ denote the upper half plane.   A {\it Fuchsian group}
is a discrete subgroup of Aut$(\mathbb{H})$.   A torsion-free Fuchsian
group is called a  {\it Fuchsian surface group}.  Fuchsian groups $\Gamma$ with
compact orbit space $\mathbb{H}/\Gamma$ have the following presentation: 
\begin{equation}\label{E:FuchPres}\langle a_1, b_1, a_2, b_2, \ldots, a_{g_0},\  b_{g_0}, \ c_1, \ldots, c_r \  | \
c_1^{m_1}, c_2^{m_2}, \ldots, c_r^{m_r}, \  \prod\limits_{i=1}^{g_0} [a_i,
b_i], \ \prod_{j=1}^{r} c_j \rangle\end{equation} 
where the generators are in Aut$(\mathbb{H})$  and the $m_i$ are integers
satisfying $2 \leq m_1 \leq m_2 \leq \ldots \leq m_r.$ For any Fuchsian
group $\Gamma$, the numbers $g_0, r$, and the $m_i$ are uniquely
determined. The tuple $(g_0; m_1, m_2, \ldots, m_r)$ is called the {\it
  signature} of $\Gamma$ and $g_0$ is called the {\it orbit genus}
(i.e.\ the genus of the orbit space $\mathbb{H}/\Gamma$).

Fuchsian groups may be used to classify automorphism groups of Riemann
surfaces.  A finite group $G$ is the automorphism
group of a compact Riemann surface of genus $g$ at least 2 precisely when
$G$ is isomorphic to a quotient $\Gamma/K$ where $\Gamma$ is a Fuchsian
group  with a compact orbit space, and $K$ is normal in $\Gamma$,
alternatively when there is a surjective homomorphism $\Phi: \Gamma \to G$. The $c_j$ values in the presentation of $\Gamma$ correspond to the
ramification of the covering $X \to X/G$.  We refer to the tuple of
elements in $G$ \begin{equation}\label{E:genvec}\left(\Phi(a_1), \Phi(b_1), \ldots, \Phi(a_{g_0}), \Phi(b_{g_0}), \Phi(c_1),
\ldots, \Phi(c_r)\right)\end{equation}  as the
{\it generating vector} for this action (i.e.\ the image of the generators in
\eqref{E:FuchPres} under $\Phi$).  

An inequality of Poincar\'{e} provides combinatorial restrictions on the
possible $m_i$ values depending on $g$ (and thus a restriction on possible
signatures).   Combining this with the Riemann-Hurwitz formula and several
group theoretical results (for instance a generalization of a result of
Harvey \cite[Theorem 4]{harvey} for abelian groups \cite[Theorem 9.1]{breuer}),
Breuer's algorithm first generates a list of all possible signatures for
Fuchsian groups $\Gamma$ for a  given genus $g$ and given order $n$ of the
automorphism group. To speed up later searches, the algorithm also identifies situations where the possible groups associated to a fixed signature must be non-solvable or must be abelian.  

Next the algorithm searches the small group database in GAP and uses  group theoretic results to construct a list of groups $G$ of order $n$ which could have one of the determined  admissible signatures for that $n$.  Each $m_i$ in the signature represents the order of the elements in some conjugacy class in the group. If a group of order $n$ does not have the proper conjugacy class structure, it is removed from the list of potential automorphism groups.   Restrictions on the possible groups associated to a signature (such as only abelian groups or only non-solvable groups) are used at this point to limit the search to groups that satisfy these restrictions.

Finally, the algorithm determines which possible groups $G$ satisfy the
condition that there is a surjective morphism $\Phi: \Gamma \to G$.  This
step in the algorithm utilizes several different group theoretic results
concerning the structure of conjugacy classes.  A result of Scott
\cite[Theorem 1]{scott} gives a sufficient condition on the irreducible characters of
a group $G$ to show there is not a surjective homomorphism $\Phi: \Gamma
\to G.$ Conversely, abelian invariants are used to give a condition on the
existence of surjective images by determining if there is a surjective
morphism $\Phi': \Gamma/\Gamma' \to G/G'$ where $\Gamma'$ and $G'$ are the
commutator subgroups of $\Gamma$ and $G$, respectively.

For each genus $g$ from 2 to 48, the output of Breuer's code consists of
the list of automorphism groups along with their signatures.  During the
execution of the algorithm,  at the step where the testing of possible
groups occurs, the specific generating vector are computed, although
Breuer's original code does not output these values.

\section{The New Code}\label{S:code}

In this section, we focus on the new pieces added to Breuer's code in order
to
produce the generating vectors.  Of course, in
the transition from GAP to Magma, several small additional modifications
were made, but we kept most of Breuer's program intact.  In particular, we
maintain his names of functions and variables  when possible.

The starting piece of the original code is a function Breuer calls {\it RepresentativesEpimorphisms}.
Given a group and signature, it will output records  containing the
signature, the list of conjugacy classes from which the branching data
arises, the group, and the image of the generators of the corresponding
Fuchsian group as in \eqref{E:genvec}. Breuer's algorithm finds all examples of such
generating vectors up to {\it simultaneous conjugation}:  if
\eqref{E:genvec} is one generating vector, then so is
$$\left(\Phi(a_1)^h, \Phi(b_1)^h, \ldots, \Phi(a_{g_0})^h, \Phi(b_{g_0})^h, \Phi(c_1)^h,
\ldots, \Phi(c_r)^h\right)  $$
for all $h \in G$, where by $\Phi(x)^h$ we mean conjugation of $\Phi(x)$ by
$h$. All of these generating vectors will give equivalent actions, and so
we only record one example for each equivalence class.  

In
the process of translating the code  to Magma, it was necessary to
convert all groups to permutation groups in order to use the command {\it
  DoubleCosetRepresentatives}. As a result, we have added a function called {\it ConvertToPerm}, based on an
algorithm found on Magma's web page, which will convert a group of type
\verb+GrpPC+ (these are finite solvable groups given by a power-conjugate presentation)  to
a permutation group. All groups from the SmallGroup database in Magma are
given as either permutation groups (if insoluble) or type \verb+GrpPC+ (if soluble). It is important to note here that all outputs of
the modified code in Magma will be for a permutation group isomorphic to
the original group.  This includes the order in which the conjugacy classes are listed,
which can change depending on the way the group is represented in
Magma. Also, the code returns an error if the group entered is not already
a permutation group or else a group of type \verb+GrpPC+.

\vspace{.1 in}

  Instead of running Breuer's full program from scratch, we used his original data and entered that into the Magma version of
the modified function {\it RepresentativesEpimorphisms}. To do this, we
added a new function called {\it AddMonodromy} which inputs a
genus $g$ and reads the Breuer data files for that particular genus.  Then
it calls the {\it RepresentativesEpimorphisms} function for each group and
signature pairing and outputs the full list of groups, signatures, and
the generators of the corresponding Fuchsian group. 

\vspace{.1 in}

We did make a few small modifications to Breuer's data before
running the new code on this data.  Because of
the limitations of GAP's database fifteen years ago, the program could not
compute a few special cases where the automorphism group had size larger
than $1000$ or size $512$ or $768$. For the few groups that exceeded the
capacity of the small group database at the time, Breuer computed
possible groups by hand using the higher genus algorithm.  He  included a special file of large group data
for the rest of the program to access in these special cases. 

Today the database of small groups in GAP and Magma includes all groups of
order up to $2000$. We have now run our Magma version of Breuer's full code on those special cases
(with two exceptions mentioned below)  for genus up to $48$ and added the group information
for those entries (i.e.\ the group identifier number) to the data files
Breuer originally
generated.   The two cases which still exceed the current database
of known groups  are groups of order $2160$ in genus 46 with
signature $(0; 2, 3, 8)$. They are identified in the Atlas of Finite Group
Representations \cite{Atlas}  as  $3.A_6.2$, which is a  maximal
subgroup of the Janko group $J_2$,  and  the central product of $S_3$
and $A_6.2$. These two groups are not included with the other genus $46$
data,  but generators and relations for the groups, as
well as the generating vectors for these actions are provided in a separate
file called {\it g46$\_$2160}. 

  Note that in Breuer's original work
there was one missing entry for genus $33$ of size $1536$ which is
mentioned in the  errata for his book \cite{breuererrata}. Additionally we have found in genus
41 that there are two groups of order $768$, not one,
with signature $(0;2,3,16)$. Both groups $( 768, 1085329 )$ and $( 768,
1085344 )$, listed with their group identification numbers from the
databases in Magma or GAP, have elements satisfying the ramification
type. We added these entries to Breuer's original data.

Finally, Breuer's representation of groups of order $256$ was not in the
standard  group identification form from the database of small order groups
in Magma or GAP. We have rerun his program
on groups of that order, identified the admissible groups as such, and manually added that
information to his original data files.

\section{Program Files and Data Files}\label{S:data}

Below is a list of the files and data sets packaged with this paper. Some
details about the functionality may be found as comments accompanying the
software.  All files may be found at:

\verb+http://www.math.grinnell.edu/+$\sim$\verb+paulhusj/monodromy.html+.

\subsubsection*{genvectors.m}

 This file contains the pertinent Magma code used to
generate the data sets described below.  Most of the algorithms and  functions in this file are 
translations of Breuer's program to Magma. We do not include all of
Breuer's functions, rather we focus only on those which compute the
generating vectors, assuming that the group and signature are already known. Comments
which begin with a ``$\#$'' symbol are Breuer's comments from the original
code.

\subsubsection*{searchroutines.m}

This file contains programs written by the author of this paper to allow for searching of
the data files inside of Magma.  We describe the details of the functions
in this file in Section \ref{S:other}.

\subsubsection*{groupsignaturedata}

  This file gives a list of each group which appears as
a subgroup of an automorphism group of a curve of genus up to $48$ without
the generating vectors. The file does not add any new data from Breuer's original data files.  We include it here merely for
completeness. It is intended for anyone who would be interested in searching the
database and it was set up with the intention that {\it grep} commands would
work easily on it.  Additionally this data file will be  necessary for
our future project of
extending the algorithm to higher genus, which relies on recursively searching the
lower genera data.

 Breuer's original data does contain a similar file but we have added the additional data,
 as described in Section \ref{S:code}. Each row in this file consists of an entry given in the form:

\begin{center}
\verb+[*genus, order of group, signature, group identification number *]+.  
\end{center}

For example, the Hurwitz curve of genus $7$ has automorphism group
PSL$(2,8)$ with signature $(0;2,3,7)$ and is represented  in the file by the following line:
\begin{center}
\verb+[*7, 504, [0,2,3,7], ( 504,156 ) *]+.  
\end{center}

\vspace{.2 in}
  
The rest of the files house the data. Each folder described below
contains a distinct file for each genus.  The first folder contains
only those examples where the size of the group $G$ is larger that $4(g-1)$,  and should
be  sufficient for many researchers. The second folder contains all the
data for genus up to $20$.   The third file is the same data for
GAP instead of Magma.

\subsubsection*{GenVectMagma-LargeGps}

These files contain the examples of ``large'' groups, those of size greater than $4(g-1)$
for a given genus $g$. This condition on the size of the group ensures that
the orbit genus $g_0$ is zero and that there are at most four branch points
\cite[Lemma 3.18]{breuer}. For anyone who is only considering groups of this
size, these are the files to use as they are substantially smaller than the
full data. (As mentioned above, for genus $46$ there are two examples of
groups of order $2160$ acting on a curve of this genus.  Since these groups
are not in the SmallGroup database, they are treated in a separate file
called {\it g46$\_$2160}.)

The entries in the file are separated with a \verb+*+ symbol on its own
row.  Each entry consists of a number of rows which give us the group
identification number,  the signature, the numbers corresponding to the  conjugacy classes from which the
generating vectors come, and the elements of the generating vectors. Again, the group is considered as a
permutation group and so the generating vector is given as a  list of permutations,
  and the ordering of the conjugacy class is for the group as a
permutation group. 

In the case where there is no ramification in the cover, the line(s) for
the permutations representing the generating vectors are replaced by a ``[ \ ]''.

Returning to the example above, the Hurwitz curve of genus $7$ is represented by the
following lines:

\begin{verbatim}
     (504,156)
     [ 0, 2, 3, 7 ]
     [ 2, 3, 4 ]
     1 6 4 3 9 2 8 7 5
     4 5 8 9 6 2 3 7 1
     5 2 8 1 6 9 7 4 3
     *
     (504,156)
     [ 0, 2, 3, 7 ]
     [ 2, 3, 5 ]
     1 6 4 3 9 2 8 7 5
     4 8 9 6 3 1 2 7 5
     2 8 9 1 5 3 7 6 4
     *
     (504,156)
     [ 0, 2, 3, 7 ]
     [ 2, 3, 6 ]
     1 6 4 3 9 2 8 7 5
     6 4 2 3 1 5 8 9 7
     9 4 3 6 2 1 5 8 7
     *
\end{verbatim}

\noindent In this example, the first entry lists the generating vector as the
three elements of the permutation group isomorphic to PSL$(2,8)$
written as  $(2\ 6)(3\ 4)(5\
9)(7\ 8)$,
$(1\ 4\ 9)(2\ 5\ 6)(3\ 8\ 7)$, and $(1\ 5 \ 6 \ 9\ 3\ 8 \ 4)$. Notice that
in this case the three sets of generating vectors above (representing all
possible generating vectors up to simultaneous conjugation) actually generate
equivalent actions, but the generating vectors come from different
conjugacy classes.  It is not part of Breuer's code to determine actions
up to equivalence, so  we also list all possible actions found through his program.   It is future work of the author
 to determine which entries are the same up to conformal and topological
 equivalence. 

\subsubsection*{GenVectMagmaToGenus20} 

\vspace{.1 in}

The folder {\it GenVectMagmaToGenus20} contains files, one for each genus from $2$ to
$20$ with the generating vectors up to simultaneous conjugation included. For genus greater than $20$, recording all the possible generating vectors  for a given group and signature
creates a huge file.  For example, given the abelian group
$\Z/2\Z \times \cdots \times \Z/2\Z$ and a signature of the form $[0;2, 2,
\ldots, 2]$ with a large number of $2$'s (say $r$ of them),  the program
finds  many
possible lists of $r$ elements of $G$ such that $g_1 \cdots g_r=1_G$ and
the $g_i$ generate the group.  Each of these would correspond to a distinct
entry of length $r+3$ in our files.  As $r$ grows large, this produces huge
data files.

 The full data for genus beyond $20$ is available upon request from the author if
a complete list of generating vectors is necessary, or
the command {\it RepresentativesEpimorphisms},  described in Section
\ref{S:other}, may be used to determine generating vectors for specific cases.

\subsubsection*{GenVectGAP}

\vspace{.1 in}

These files are the GAP versions of the generating vector data.  Unlike in the
Magma version, we have rerun Breuer's full program instead of just
computing the generating vectors from his data. Since the program
we used to generate this data is identical to Breuer's
original program, except to add an extra command to print the generating vectors, we
have not included the code for this program. The output files contain rows with the following information:

\verb+[ group ][ signature ][ conjugacy classes ][ generating vector ]+.

\noindent The generating vectors are given as elements of the particular 
group. Here is the output for the Hurwitz curve in genus $7$:

\begin{center}

\SMALL \verb+[ 504, 156 ][ 0, 2, 3, 7 ][ 5, 6, 2 ][ (2,3)(4,6)(5,8)(7,9), (1,2,9)(3,4,6)(5,8,7), (1,7,5,9,3,4,2) ]+
\verb+[ 504, 156 ][ 0, 2, 3, 7 ][ 5, 6, 3 ][ (2,3)(4,6)(5,8)(7,9), (1,2,4)(3,5,8)(6,9,7), (1,6,9,4,3,5,2) ]+
\verb+[ 504, 156 ][ 0, 2, 3, 7 ][ 5, 6, 4 ][ (2,3)(4,6)(5,8)(7,9), (1,2,5)(3,9,7)(4,6,8), (1,8,4,5,3,9,2) ]+.
\end{center}

There is also a software package in GAP  called ``MapClass'', written by
James, Magaard, Shpectorov, and  V\"{o}lklein, which, among other computations,
will find the generating vectors given a group and list of conjugacy
classes corresponding to a signature \cite{braid}.

\section{Additional Commands}\label{S:other}

For any curve (regardless of genus), if the signature and automorphism
group is already known,  the command {\it RepresentativesEpimorphisms(sign,
group)}, found in the file {\it genvectors.m}, will
compute the corresponding generating vectors. If the group is not given as a
permutation group or a group of type \verb+GrpPC+ in Magma, then an error is
returned. The output of this command is a list of records. Each record
contains four fields.  They are as follows:

\verb+signature+  the signature,

\verb+Con+  the number of the conjugacy classes where the generating
vectors come from (in Magma as group of type \verb+GrpPerm+),
 
\verb+Gro+   the group as a permutation group, and
   
\verb+genimages+  the generators of the automorphism group corresponding to
the image of the generators of the Fuchsian group as in \eqref{E:FuchPres} under the
map $\Phi: \Gamma \to G$; the first $2g_0$
elements correspond to the image of the $a_i$ and $b_i$ values in \eqref{E:FuchPres}, and
the rest of the entries are the images of the branching data, all given as
a list of permutation elements in $G$. \\

For instance, Hurwitz curves have signature $[0;2,3,7]$ and there is
an infinite family of these curves with automorphism group PSL$(2,q)$ where
$q$ is a prime congruent to $\pm 1 \bmod 7$ \cite[Theorem 8]{macbeath}. There is a
curve of genus $146$ with automorphism group PSL$(2,29)$ and corresponding signature $[0;
2, 3, 7]$.  If we run the
command

\verb+RepresentativesEpimorphims([0,2,3,7], PSL(2,29))+

\noindent the program will output a list of six records, two for each of
three different lists of conjugacy classes, depending on which of the three
conjugacy classes of elements of order $7$ are considered.

\vspace{.1 in}

The file {\it searchroutines.m} contains a program which can be called in
Magma to search the data, as well as some sample programs to demonstrate how one
might use this program.  The key routine is {\em ReadData(filename,
function)}  where {\it filename} should be the corresponding data file
(from those described in the previous section) for
a particular genus, and  {\it function} is a function which takes three inputs: a
group, a signature, and the generating vectors, and returns a boolean.
{\it ReadData} will read each entry of the data file of a particular genus and
then run {\it function} on each entry. If the boolean function returns a
\verb+true+ value, the data from that particular entry (the group, signature, and
generating vector) will be saved in a list, and {\it ReadData} will return the
full list after it tests each entry of the data for a given genus.

Several sample functions which use the {\it ReadData} function are also
provided in this file, including:  

\verb+FindGroup(filename,gpsize,gpnum)+   returns all groups matching the
group 

defined as \verb+SmallGroup(gpsize,gpnum)+,

\verb+FindSignature(filename, signature)+  returns all entries with
signature 

matching the second input value, and

\verb+LoadAllData(filename)+   reads all data from the file,  stored in
        a list of Magma 

readable information. This will not be tenable for high 
          genus, as the data is too 

large.

\vspace{.1 in}

Finally, there are two samples of boolean functions which may be used to input into
{\it ReadData}.  Both samples will determine if an entry of the database
represents a nilpotent group.  The first sample will save all of the
nilpotent examples as a
list to Magma, while the second will write the data to a specified output
file.  

A warning: the group will be read into ReadData as \verb+SmallGroup(a,b)+, but
to actually apply the elements of a generating vector, you will need to
convert the group to a permutation group first (code is provided to make
that conversion in this file--see the function {\it ConvertToPerm(G)} in
{\it genvectors.m} to convert from a group of type \verb+GrpPC+ to a
permutation group).

\section{Acknowledgments}

The author wishes to thank Michael Zieve who originally connected her with
Breuer's program and  discussed the progress of this modified program several times. Also, thanks to Syed Lavasani who beta
tested an early version of the low genus GAP data.  Additionally, the author is 
grateful for a series of
helpful conversations with  Anita Rojas, David Swinarski, and Aaron Wooton.

Finally, thanks to Thomas Breuer for permission to make publicly available
the modifications of his code, and the corresponding data.

\end{document}